\documentclass[12pt,leqno]{amsart}
\usepackage{amsmath,amssymb,amscd,latexsym,amsthm,mathrsfs}
\usepackage{color}
\usepackage[unicode]{hyperref}
\usepackage{hypbmsec,longtable}
\ifx\pdfoutput\undefined\else
\hypersetup{
colorlinks=true,
linkcolor=blue,
citecolor=blue,
urlcolor=blue,
filecolor=blue,
bookmarksnumbered=true,
pdfstartview=FitH,
pdfhighlight=/N
}
\fi

\newtheorem{theorem}{Theorem}

\newtheorem{conjecture}[theorem]{Conjecture}

\theoremstyle{remark}

\newtheorem*{acknowledgements}{Acknowledgements}
\renewcommand{\d}{{\mathrm d}}
\renewcommand{\Im}{\operatorname{Im}}
\renewcommand{\Re}{\operatorname{Re}}

\begin{document}

\hypersetup{pdfauthor={Mathew Rogers, Wadim Zudilin},%
pdftitle={Modular equations and lattice sums}}

\title{Modular equations and lattice sums}

\author{Mathew Rogers}
\address{Department of Mathematics, University of Illinois, Urbana, IL 61801, USA}
\email{mathewrogers@gmail.com}

\author{Boonrod Yuttanan}
\address{Department of Mathematics, University of Illinois, Urbana, IL 61801, USA}
\email{byuttan2@illinois.edu}

\thanks{The first author is supported by National Science Foundation award DMS-0803107.}

\date{July 28, 2011}


\begin{abstract}
We highlight modular equations discovered by Somos and Ramanujan, and use them to
prove new relations between lattice sums and hypergeometric
functions. We also discuss progress towards solving Boyd's Mahler
measure conjectures, and conjecture a new formula for $L(E,2)$ of conductor $17$ elliptic curves.
\end{abstract}

\maketitle

\section{Introduction}
\label{s-intro}

Modular equations appear in a variety of number-theoretic contexts.
Their connection to formulas for $1/\pi$ \cite{Ra}, Ramanujan
constants such as $e^{\pi\sqrt{163}}$ \cite{We}, and elliptic curve
cryptography is well established.  In the classical theory of
modular forms, an $n$th degree modular equation is an algebraic
relation between $j(\tau)$ and $j(n\tau)$, where $j(\tau)$ is the
$j$-invariant.  For our purposes a modular equation is simply a
non-trivial algebraic relation between theta (or eta) functions. In
this paper we use modular equations to study four-dimensional lattice
sums. The lattice sums are interesting because they arise in the
study of Mahler measures of elliptic curves.

There are a large number of hypothetical relations between special values of $L$-series of elliptic curves, and Mahler measures of two-variable polynomials.
The Mahler measures $m(\alpha)$, $n(\alpha)$, and
$g(\alpha)$ are defined by
\begin{equation}\label{Mahler measure definitions}
\begin{split}
m(\alpha):=&\int_{0}^{1}\int_{0}^{1}\log\left|y+y^{-1}+z+z^{-1}+\alpha\right|\d\theta_1
\d\theta_2,\\
n(\alpha):=&\int_{0}^{1}\int_{0}^{1}\log\left|y^3+z^3+1-\alpha y
z\right|\d\theta_1\d\theta_2,\\
g(\alpha):=&\int_{0}^{1}\int_{0}^{1}\log\left|(y+1)(z+1)(y+z)-\alpha
y z\right|\d\theta_1\d\theta_2,
\end{split}
\end{equation}
where $y=e^{2\pi i\theta_1}$, and $z=e^{2\pi i\theta_2}$.  Boyd conjectured that for all integral values of
$k\not=4$ \cite{Bo1}:
\begin{equation*}
m(k)\stackrel{?}{=}\frac{q}{\pi^2} L(E,2),
\end{equation*}
where $E$ is an elliptic curve, $q$ is rational, and both $E$ and $q$ depend on $k$.  He also
discovered large numbers of formulas involving $g(\alpha)$ and $n(\alpha)$.  In cases where $E$ has a small
conductor, it is frequently possible to express $L(E,2)$ in terms of four-dimensional lattice sums.  Thus many of Boyd's identities can be regarded as series
acceleration formulas.  The main goal of this paper is to prove new formulas for the lattice sum $F(b,c)$,
defined in \eqref{F(b,c) definition}.  So far there are at least $18$ instances where $F(b,c)$ is known (or conjectured) to reduce to integrals of elementary functions.
The modular equations of Somos and Ramanujan are the main tools in our analysis.
\section{Eta function product identities}

Somos discovered thousands of new modular equations by searching for linear
relations between products of Dedekind eta functions.  Somos refers to these formulas as ``eta function product identities".  The existence
of eta function product identities can be established by using fact
that certain modular parameters (such as $j(\tau)$) equal
rational expressions involving eta functions.
By clearing denominators it is possible to rewrite
classical modular equations as eta function product identities.  A major surprise of Somos's experimental approach,
is that it turned up a large number of unexpectedly simple
identities. In order to give an example, first consider the eta
function with respect to $q$:
\begin{equation*}
\eta(q)=q^{1/24}\prod_{n=1}^{\infty}(1-q^n)=\sum_{n=-\infty}^{\infty}(-1)^n
q^{(6n+1)^2/24},
\end{equation*}
and adopt the short hand notation
\begin{equation*}
e_j=\eta(q^j).
\end{equation*}
The following formula is the smallest eta function product identity
in Somos's list \cite{so2}:
\begin{equation}\label{Somos three term}
e_{2}e_{6}e_{10}e_{30}=e_{1}e_{12}e_{15}e_{20}+e_{3}e_{4}e_{5}e_{60}.
\end{equation}
Notice that all three monomials are products of four eta functions,
and are essentially weight-two modular forms. There are no known
identities between eta products of weight less than two, and
\eqref{Somos three term} appears to be the only three-term linear
relation between products of four eta functions. Many additional
identities are known if the number of terms is allowed to
increase, or if eta products of higher weight are considered. For additional examples see
formulas \eqref{modular eq. degree 18}, \eqref{modular equation for F(4,7,7,28)}, \eqref{deg 25 modular equation}, \eqref{eta prod deg 45}, and \eqref{Somos's second smallest identity}.

Identities such as \eqref{Somos three term} can be proved almost
effortlessly with the theory of modular forms. A typical proof
involves checking that the first few Fourier coefficients of a
presumed identity vanish. Sturm's Theorem furnishes an upper bound on the
number of coefficients that need to be examined \cite{Ono2}.
We note that it is often possible, but usually more difficult, to
prove such identities via elementary $q$-series methods. Ramanujan
was a master $q$-series manipulator, and his notebooks are filled
with various modular equations and their corollaries.  We conclude
this section by providing a Ramanujan-style proof of \eqref{Somos
three term}. The main news is that \eqref{Somos three
term} can be derived from modular equations known to Ramanujan.

\begin{theorem}\label{Theorem about Somos's smallest identity} The following formula is true:
\begin{equation}\label{Somos's smallest identity}
e_{2}e_{6}e_{10}e_{30}=e_{1}e_{12}e_{15}e_{20}+e_{3}e_{4}e_{5}e_{60}.
\end{equation}
\end{theorem}
\begin{proof} Before proving \eqref{Somos's smallest identity} we need to define a small amount of
notation.  Let us denote the usual theta
functions by
\begin{align}\label{theta definition}
\varphi(q):=&\sum_{n=-\infty}^{\infty}q^{n^2},
&&\psi(q):=\sum_{n=0}^{\infty}q^{n(n+1)/2}.
\end{align}
Furthermore define $u_j$ and $z_j$ by
\begin{align*}
u_j:=&1-\frac{\varphi^4(-q^j)}{\varphi^4(q^j)},&&z_j:=\varphi^2(q^j).
\end{align*}
Notice that Ramanujan uses a slightly different notation \cite{Be3}.  
He typically sets $\alpha=u_1$, and he says that ``$\beta$ has degree $j$ over $\alpha$" when $\beta=u_j$.
Certain values of the eta function can be
expressed in terms of $u_1$ and $z_1$ \cite[p.~124]{Be3}.  We
have
\begin{align}
\eta(q)=&2^{-1/6}u_1^{1/24}(1-u_1)^{1/6}\sqrt{z_1},\label{etainversion1}\\
\eta(q^2)=&2^{-1/3}\left\{u_1(1-u_1)\right\}^{1/12}\sqrt{z_1},\label{etainversion2}\\
\eta(q^4)=&2^{-2/3}u_1^{1/6}(1-u_1)^{1/24}\sqrt{z_1}.\label{etainversion3}
\end{align}
Now we prove \eqref{Somos's smallest identity}.  By
\eqref{etainversion2} the left-hand side of the identity becomes
\begin{equation*}
e_2e_6e_{10}e_{30}=2^{-4/3}\{u_1 u_3 u_5 u_{15}(1-u_1)(1-u_{3})(1-u_{5})
(1-u_{15})\}^{1/12}\sqrt{z_{1}z_{3}z_{5}z_{15}}.
\end{equation*}
By \eqref{etainversion1} and \eqref{etainversion3}, the right-hand
side of the identity becomes
\begin{equation*}
\begin{split}
e_1 e_{12}e_{15}e_{20}&+e_3 e_4 e_5 e_{60}\\
=&2^{-5/3}\left(\left\{
u_3 u_5(1-u_1)(1-u_{15})\right\}^{1/6}\left\{
u_1 u_{15}(1-u_{3})(1-u_{5})\right\} ^{1/24}\right.\\
&+\left.\left\{ u_1 u_{15}(1-u_3)(1-u_5)\right\}
^{1/6}\left\{ u_3 u_5(1-u_1)(1-u_{15})\right\}
^{1/24}\right)\sqrt{z_{1}z_{3}z_{5}z_{15}}.
\end{split}
\end{equation*}
Combining the last two formulas shows that \eqref{Somos's smallest
identity} is equivalent to
\begin{equation}\label{eq thm 2.1 midway}
\begin{split}
2^{1/3}&\left\{u_1 u_3 u_5 u_{15}(1-u_1)(1-u_3)(1-u_{5})(1-u_{15})\right\}^{1/24}\\
&=\left\{ u_3 u_5(1-u_1)(1-u_{15})\right\}
^{1/8}+\left\{ u_1 u_{15}(1-u_3)(1-u_5)\right\}
^{1/8}.
\end{split}
\end{equation}
It is sufficient to show that \eqref{eq thm 2.1 midway} can be
deduced from Ramanujan's modular equations.  

The first modular equation we require can be recovered by
multiplying entries 11.1 and 11.2 in \cite[p.~383]{Be3}:
\begin{equation*}
\left((u_1 u_{15})^{1/8}+\left\{
(1-u_1)(1-u_{15})\right\} ^{1/8}\right)\left(
(u_3 u_5)^{1/8}+\left\{ (1-u_3)(1-u_5)\right\}
^{1/8}\right)=1.
\end{equation*}
Rearranging yields an
identity for the right-hand side of \eqref{eq thm 2.1 midway}:
\begin{equation}\label{eq:2-1-intermed}
\begin{split}
&\{u_3 u_5(1-u_1)(1-u_{15})\}^{1/8}+\{u_1 u_{15}(1-u_3)(1-u_5)\}^{1/8}\\
&\qquad=
1-\{u_1 u_3 u_5 u_{15}\}^{1/8}-\{(1-u_1)(1-u_3)(1-u_5)(1-u_{15})\}^{1/8}.
\end{split}
\end{equation}
By entry 11.14 in \cite[p.~385]{Be3}, it is clear that
\begin{equation}\label{eq:2-1-intermed2}
\begin{split}
&1-\{u_1 u_3 u_5 u_{15}\}^{1/8}-\{(1-u_1)(1-u_3)(1-u_5)(1-u_{15})\}^{1/8}\\
&\qquad=
2^{1/3}\{u_1 u_3 u_5 u_{15}(1-u_1)(1-u_3)(1-u_5)(1-u_{15})\}^{1/24}.
\end{split}
\end{equation}
The theorem follows from combining \eqref{eq:2-1-intermed} and
\eqref{eq:2-1-intermed2} to recover \eqref{eq thm 2.1 midway}.
\end{proof}

It appears to be very difficult to explain \textit{why} identities such as
\eqref{Somos's smallest identity} exist.  Our initial motivation for
constructing an elementary proof of \eqref{Somos's smallest
identity}, was to find a method for generating more
identities. It would be interesting if a systematic
method for generating weight-two eta product identities could be discovered.  This is a surprisingly important question for studying lattice
sums and the Beilinson conjectures.

\section{Lattice Sums}\label{Section:Lattice sums}

In this section we investigate four-dimensional lattice
sums.  Many of these sums appear in the study of Mahler measures of elliptic curves.
Let us define
\begin{equation*}
\begin{split}
F(a,b,c,d):=&(a+b+c+d)^2\\
&\times\sum_{n_i=-\infty}^{\infty}\frac{(-1)^{n_1+n_2+n_3+n_4}}{\left(a(6n_1+1)^2+b(6n_2+1)^2+c(6n_3+1)^2+d
(6n_4+1)^2\right)^2}.
\end{split}
\end{equation*}
The four-dimensional series is not absolutely convergent, so it
is necessary to employ
summation by cubes \cite{Bor}.  Notice that Euler's pentagonal number theorem
can be used to represent $F(a,b,c,d)$ as a useful
integral
\begin{equation}\label{F(a,b,c,d) integral formula}
F(a,b,c,d)=-\frac{(a+b+c+d)^2}{24^2}\int_{0}^{1}\eta(q^a)\eta(q^b)\eta(q^c)\eta(q^{d})
\log q\frac{\d q}{q}.
\end{equation}
We also use the shorthand notation
\begin{equation}\label{F(b,c) definition}
F(b,c):=F(1,b,c,b c),
\end{equation}
since we are primarily interested in cases where $a=1$, $d=b c$, and
$b$ and $c$ are rational.

The interplay between values of $F(b,c)$, Boyd's Mahler measure
conjectures, and the Beilinson conjectures is outlined in
\cite{Rgsubmit}.  If $(b,c)\in\mathbb{N}^2$ and $(1+b)(1+c)$ divides
$24$, then $F(b,c)=L(E,2)$ for an elliptic curve $E$. Formulas are now rigorously proved relating each of those eight
cases to Mahler measures such as $m(\alpha)$ \cite{Zu}.  Because Mahler measures often reduce to generalized hypergeometric
functions, many of Boyd's identities can be regarded as series transformations  \cite{RV}, \cite{LR}. It is known that
\begin{align*}
m(\alpha)=&\Re\left[\log(\alpha)-\frac{2}{\alpha^2}{_4F_3}\left(\substack{\frac{3}{2},\frac{3}{2},1,1\\
2,2,2};\frac{16}{\alpha^2}\right)\right],\text{if $\alpha\not=0$,}\\
n(\alpha)=&\Re\left[\log(\alpha)-\frac{2}{\alpha^3}{_4F_3}\left(\substack{\frac{4}{3},\frac{5}{3},1,1\\
2,2,2};\frac{27}{\alpha^3}\right)\right],\text{ if $|\alpha|$ is
sufficiently large,}\\
3g(\alpha)=&n\left(\frac{\alpha+4}{\alpha^{2/3}}\right)+4n\left(\frac{\alpha-2}{\alpha^{1/3}}\right),\text{ if $|\alpha|$ is
sufficiently large.}
\end{align*}
The function $m(\alpha)$ also reduces to a $_3F_2$ function
if $\alpha\in\mathbb{R}$ \cite{KO}, \cite{Rgsubmit}.  The first author and Zudilin \cite{RZ1}
recently proved that
\begin{equation}\label{deninger formula}
F(3,5)=\frac{4\pi^2}{15}m(1)=\frac{\pi^2}{15}{_3F_2}\left(\substack{\frac{1}{2},\frac{1}{2},\frac{1}{2}\\1,\frac{3}{2}}\bigg|\frac{1}{16}\right).
\end{equation}
Equation \eqref{deninger formula} is equivalent to a formula that
was conjectured by Deninger \cite{De}, and which helped motivate
Boyd's seminal paper \cite{Bo1}.  One of the main results in
\cite{Rgsubmit}, is that it is also possible to find formulas for
values such as $F(1,4)$ and $F(2,2)$.  These cases are probably not
related to elliptic curve $L$-values.  As a result it was
hypothesized that it should be possible to ``sum up" $F(b,c)$ for
arbitrary values of $b$ and $c$.
\subsection{Lacunary cases}
It is typically very difficult to prove formulas such as
\eqref{deninger formula}. The proof of \eqref{deninger formula} is a
$q$-series proof which utilizes the integral representation
\eqref{F(a,b,c,d) integral formula}.  In general the difficulty of
dealing with a lattice sum depends on whether it is
\textit{lacunary} or \textit{non-lacunary}.  The lacunary values can
be reduced to two-dimensional sums, which are (almost) always easier
to deal with than four-dimensional sums. It is often difficult to
determine whether or not a particular sum is lacunary. Cases such as
$F(1,1)$, $F(1,2)$, and $F(1,3)$ can be equated to $L$-values of CM
elliptic curves, and their lacunarity follows from the CM
hypothesis. Values such as $F(1,4)$ and $F(2,2)$ have no arithmetic
interpretation, however they easily reduce to two-dimensional sums
via classical theta series results. The usual method for detecting
lacunarity, is to expand the associated cusp form in an infinite
series.  If one writes
\begin{equation*}
\eta(q^a)\eta(q^b)\eta(q^c)\eta(q^d)=q^{(a+b+c+d)/24}\left(a_0+a_1
q+a_2 q^2+\dots\right),
\end{equation*} then the non-vanishing $a_i$'s should have
zero-density. For the cases discussed herein it is usually necessary
to compute thousands of coefficients to observe lacunarity.
Additional non-obvious lacunary values include $F(2,9)$ and
$F(4,7,7,28)$.  It is necessary to employ eta function product
identities to deal with these last two cases. By a result of
Ramanujan \cite[p.~210, Entry 56]{Be4}, we have
\begin{equation}\label{modular eq. degree 18}
3 e_1 e_2 e_9 e_{18}=-e_1^2 e_2^2+e_1^3 \frac{e_{18}^2}{e_9}+e_2^3
\frac{e_9^2}{e_{18}}.
\end{equation}
Substituting classical theta expansions for $e_1^3$, $e_2^2/e_1$,
and $e_1^2/e_2$ \cite[pg.~114-117]{GK}, leads to
\begin{equation}\label{F(2,9) 1-2-9-18 series}
\begin{split}
3\eta(q)\eta(q^2)\eta(q^9)\eta(q^{18})=&-\sum_{\substack{n=0\\k=0}}^{\infty}
(-1)^n
(2n+1)q^{\frac{(2n+1)^2+(2k+1)^2}{8}}\\
&+\sum_{\substack{n=0\\k=0}}^{\infty}
(-1)^n (2n+1)q^{\frac{(2n+1)^2+9(2k+1)^2}{8}}\\
&+\sum_{\substack{n=0\\
k=-\infty}}^{ \infty}(-1)^{n+k}(2n+1)q^{\frac{(2n+1)^2+9(2k)^2}{4}}.
\end{split}
\end{equation}
Because $e_1 e_2 e_9 e_{18}$ is a finite linear combination of
two-dimensional theta series, it must be a lacunary eta product.
Formula \eqref{F(2,9) 1-2-9-18 series} is the main ingredient needed
to relate $F(2,9)$ to hypergeometric functions and Mahler measures.

\begin{theorem}\label{lattice sums main theoerem} Let $t=\sqrt[4]{12}$, then the following identity is true:
\begin{equation}\label{F(2,9) closed form}
\begin{split}
\frac{144}{25\pi^2}F(2,9)=&-3m\left(4i\right)+2m\left(\frac{1}{\sqrt{2}}\left(4-2t-2t^2+t^3\right)\right)\\
&+m\left(4i\left(7+4t+2t^2+t^3\right)\right).
\end{split}
\end{equation}
\end{theorem}
\begin{proof}  The most difficult portion of the calculation is
to find a two-dimensional theta series for $e_1 e_2 e_9 e_{18}$.
This task has been accomplished via an eta function product
identity.  The remaining calculations parallel those carried out in
\cite{Rgsubmit}. Integrating \eqref{F(2,9) 1-2-9-18 series} leads to
\begin{equation}\label{integrated modular equation}
\begin{split}
\frac{3}{25}F(2,9)+F(1,2)=&4\sum_{\substack{n=0\\k=0}}^{\infty}
\frac{(-1)^n (2n+1)}{\left((2n+1)^2+9(2k+1)^2\right)^2}\\
&+\sum_{\substack{n=0\\
k=-\infty}}^{\infty}\frac{(-1)^{n+k}(2n+1)}{\left((2n+1)^2+9(2
k)^2\right)^2}.
\end{split}
\end{equation}
There are two possible formulas for $F(1,2)$ \cite{RV}:
\begin{align}
F(1,2)=&\frac{\pi^2}{8}m\left(2\sqrt{2}\right)
=\frac{\pi^2}{16}m\left(4i\right)\label{F(1,2)}.
\end{align}
By the formula for $F_{(1,2)}(3)$ in \cite[Eq.~115]{Rgsubmit}, we also have
\begin{equation}\label{lattice 1}
\sum_{\substack{n=0\\
k=-\infty}}^{\infty}\frac{(-1)^{n+k}(2n+1)}{\left((2n+1)^2+9(2
k)^2\right)^2}=\frac{\pi^2}{48}m\left(4i\left(7+4t+2t^2+t^3\right)\right),
\end{equation}
where $t=\sqrt[4]{12}$. Next we evaluate the remaining term in
\eqref{integrated modular equation}.  Notice that for $x>0$
\begin{equation*}
\begin{split}
\sum_{\substack{n=0\\k=0}}^{\infty}& \frac{(-1)^n
(2n+1)}{\left((2n+1)^2+x(2k+1)^2\right)^2}\\
=&\frac{\pi^2}{16}\int_{0}^{\infty}u\left(\sum_{\substack{n=0}}^{\infty}
(-1)^n (2n+1)e^{-\pi(n+1/2)^2
u}\right)\left(\sum_{k=0}^{\infty}e^{-\pi x(k+1/2)^2 u}\right)\d u.
\end{split}
\end{equation*}
By the involution for the weight-$3/2$ theta function
\begin{equation*}
\sum_{n=0}^{\infty}(-1)^n(2n+1)e^{-\pi(n+1/2)^2
u}=\frac{1}{u^{3/2}}\sum_{n=0}^{\infty}(-1)^n
(2n+1)e^{-\pi(n+1/2)^2\frac{1}{u}},
\end{equation*}
this becomes
\begin{equation*}
\begin{split}
\sum_{\substack{n=0\\k=0}}^{\infty}& \frac{(-1)^n
(2n+1)}{\left((2n+1)^2+x(2k+1)^2\right)^2}\\
&=\frac{\pi^2}{16}\sum_{\substack{n=0\\
k=0}}^{\infty}(-1)^n (2n+1)\int_{0}^{\infty}u^{-1/2}
e^{-\pi\left((n+1/2)^2 \frac{1}{u}+x(k+1/2)^2 u \right)}\d u\\
&=\frac{\pi^2}{16\sqrt{x}}\sum_{\substack{n=0\\
k=0}}^{\infty}(-1)^n \frac{(2n+1)}{(2k+1)}e^{-\frac{\pi\sqrt{x}}{2}(2n+1)(2k+1)}\\
&=\frac{\pi^2}{16\sqrt{x}}\sum_{n=0}^{\infty}(-1)^n
(2n+1)\log\left(\frac{1+e^{-\pi\sqrt{x}(n+1/2)}}{1-e^{-\pi\sqrt{x}
(n+1/2)}}\right).
\end{split}
\end{equation*}
Applying formulas (1.6), (1.7), and (2.9) in \cite{LR}, we have
\begin{equation*}
\begin{split}
&=\frac{\pi^2}{32\sqrt{x}}\left(m\left(\frac{4}{\sqrt{\alpha_{x/4}}}\right)-m\left(\frac{4i\sqrt{1-\alpha_{x/4}}}{\sqrt{\alpha_{x/4}}}\right)\right)\\
&=\frac{\pi^2}{32\sqrt{x}}m\left(4\left(\frac{1-\sqrt{1-\alpha_{x/4}}}{1+\sqrt{1-\alpha_{x/4}}}\right)\right),
\end{split}
\end{equation*}
where $\alpha_{x}$ is the singular modulus (recall that
$\alpha_x=1-\varphi^4(-e^{-\pi\sqrt{x}})/\varphi^4(e^{-\pi\sqrt{x}})$).
The second degree modular equation shows that
\begin{equation*}
\frac{1-\sqrt{1-\alpha_{x/4}}}{1+\sqrt{1-\alpha_{x/4}}}=\sqrt{\alpha_x},
\end{equation*}
and hence we obtain
\begin{equation}\label{closed form sum}
\sum_{\substack{n=0\\k=0}}^{\infty} \frac{(-1)^n
(2n+1)}{\left((2n+1)^2+x(2k+1)^2\right)^2}=\frac{\pi^2}{32\sqrt{x}}m\left(4\sqrt{\alpha_x}\right).
\end{equation}
It is well known that $\alpha_n$ can be expressed in terms of
 class invariants if $n\in\mathbb{Z}$:
\begin{equation*}
\alpha_n=\frac{1}{2}\left(1-\sqrt{1-1/G_n^{24}}\right),
\end{equation*}
and the values of $G_n$ have been extensively tabulated, read: \cite[p.
188]{Be5}. Setting $n=9$ yields
\begin{equation*}
\begin{split}
\alpha_9=&\frac{1}{2}\left(1-\sqrt{1-\left(\frac{\sqrt{2}}{\sqrt{3}+1}\right)^{8}}\right)\\
=&\frac{1}{2}\left(1-4t+t^3\right)\\
=&\frac{\left(4-2t-2t^2+t^3\right)^2}{32},
\end{split}
\end{equation*}
where $t=\sqrt[4]{12}$.  It follows immediately that
\begin{equation}\label{lattice 2}
\sum_{\substack{n=0\\k=0}}^{\infty} \frac{(-1)^n
(2n+1)}{\left((2n+1)^2+9(2k+1)^2\right)^2}=\frac{\pi^2}{96}m\left(\frac{1}{\sqrt{2}}\left(4-2t-2t^2+t^3\right)\right).
\end{equation}
The proof of \eqref{F(2,9) closed form} can be completed by
combining \eqref{integrated modular equation}, \eqref{F(1,2)},
\eqref{lattice 1}, and \eqref{lattice 2}.
\end{proof}

In order to avoid tedious calculations we have chosen to omit the
explicit formula for $F(4,7,7,28)$ from this paper\footnote{The
formula is available in previous versions of this paper on the
Arxives}.  It suffices to say that the sum reduces to an extremely
complicated expression involving hypergeometric functions and Meijer
$G$-functions. The key modular equation used to prove lacunarity is due to Somos \cite[Entry
$q_{28,9,35}$]{so}:
\begin{equation}\label{modular equation for F(4,7,7,28)}
28 e_4 e_7^2 e_{28}=-7e_1e_7^3-\frac{e_1^5}{e_2^2}\frac{
e_{14}^2}{e_7}+8\frac{e_2^5}{e_1^2}e_{14}.
\end{equation}
By classical theta expansions \cite[pg.~114-117]{GK}, the identity can be rewritten as
\begin{equation*}
\begin{split}
28 \eta(q^4)\eta^2(q^7)\eta(q^{28})
=&-7\sum_{\substack{n=-\infty\\k=0}}^{\infty}(-1)^{n+k}(2k+1)q^{\frac{(6n+1)^2+21(2k+1)^2}{24}}\\
&-\sum_{\substack{n=-\infty\\k=0}}^{\infty}(6n+1)q^{\frac{(6n+1)^2+21(2k+1)^2}{24}}\\
&+8\sum_{n,k=-\infty}^{\infty}(-1)^{n+k}(3n+1)q^{\frac{4(3n+1)^2+7(6k+1)^2}{12}},
\end{split}
\end{equation*}
and as a result it is easy to see that $e_4 e_7^2 e_{28}$ is
lacunary.

Apart from $F(2,9)$, $F(4,7,7,28)$, and the examples discussed in \cite{Rgsubmit}, we are not aware of any additional
lacunary values of $F(a,b,c,d)$ (although they probably do exist).  It is also possible to find
two-dimensional reductions for certain linear combinations of
lattice sums, however these formulas are generally less interesting than the
previous examples. To give a single case let us briefly
consider the following modular equation \cite[Entry $x_{50,6,81}$]{so}:
\begin{equation}\label{deg 25 modular equation}
5 e_1 e_2 e_{25} e_{50}+2 e_1^2 e_2 e_{50}+2 e_1 e_2^2 e_{25}=-e_1^2
e_2^2 +e_1^3 \frac{e_{50}^2}{e_{25}}+e_2^3 \frac{e_{25}^2}{e_{50}}.
\end{equation}
All three eta quotients on the right-hand side of \eqref{deg 25
modular equation} have two-dimensional theta series expansions.  As
a result it is possible to prove that
\begin{equation}\label{F(2,25) closed form}
\begin{split}
&\frac{5}{13^2} F(2,25)+\frac{2}{9^2} F(1,1,2,50)+\frac{2}{5^2}
F\left(1,2,2,25\right)\\
&\qquad=\frac{\pi^2}{80}\left(-5m\left(4i\right)+2
m\left(4\sqrt{\alpha_{25}}\right)+m\left(4i\sqrt{\frac{1-\alpha_{25}}{\alpha_{25}}}\right)\right),
\end{split}
\end{equation}
where
$\alpha_{25}=\frac{1}{2^{13}}\left(\sqrt{5}-1\right)^8\left(\sqrt[4]{5}-1\right)^8$.
There are many additional results along the lines of \eqref{F(2,25)
closed form} which we will not discuss here.

\subsection{Non-lacunary cases}

In instances where $F(a,b,c,d)$ \textit{does not} reduce to a
two-dimensional sum, the calculations become far more difficult. The
recent proofs of formulas for $F(1,5)$, $F(2,3)$ and $F(3,5)$ are
all based upon new types of $q$-integral transformations \cite{RZ1},
\cite{RZ2}. The fundamental transformation used to prove a formula for $F(2,3)$
is
\begin{equation*}
\begin{split}
&\int_{0}^{1}q^{1/2}\psi(q)\psi(q^3)\varphi(-q^x)\varphi(-q^{3x})\log
q\frac{\d q}{q}\\
&\qquad=\frac{2\pi}{3x}\Im\int_{0}^{1}\omega
q\psi^4\left(\omega^2 q^2\right)\log\left(4
q^{3x}\frac{\psi^4(q^{12x})}{\psi^4(q^{6x})}\right)\frac{\d q}{q},
\end{split}
\end{equation*}
where $\omega=e^{2\pi i/3}$. When $x=1$ the left-hand side equals $4F(2,3)$ (to see this use
$q^{1/8}\psi(q)=\eta^2(q^2)/\eta(q)$ and
$\varphi(-q)=\eta^2(q)/\eta(q^2)$), and the right-hand side becomes an extremely complicated elementary integral.  The
most difficult portion of the calculation is to reduce the
elementary integral to hypergeometric functions. It was proved with difficulty that
\begin{equation*}
F(2,3)=\frac{\pi^2}{6}m(2)=\frac{\pi^2}{12}{_3F_2}\left(\substack{\frac{1}{2},\frac{1}{2},\frac{1}{2}\\1,\frac{3}{2}};\frac{1}{4}\right).
\end{equation*}
Boyd's numerical work was instrumental in the calculation, because it
allowed for the final formula to be anticipated in advance.

There are many cases where non-lacunary lattice sums reduce to
elementary integrals, but where the integrals are extremely
difficult to deal with.  We recently used the method from \cite{RZ1}
to find a formula for $F(1,8)$:
\begin{equation}\label{F(1,8) formula}
F(1,8)=\frac{9\pi\sqrt[4]{2}}{128}\int_{0}^{1}\frac{(1-k)^2+2\sqrt{2(k+k^3)}}{(1+k)(k+k^3)^{3/4}}
\log\left(\frac{1+2k-k^2+2\sqrt{k-k^3}}{1+k^2}\right)\d k
\end{equation}
We checked this monstrous identity to more than $100$ decimal places
by calculating $F(1,8)$ with \eqref{F(a,b,c,d) integral formula}.
We can only speculate that the integral should reduce to
something along the lines of \eqref{F(2,9) closed form}.

Occasionally eta function identities provide shortcuts for avoiding
integrals like \eqref{F(1,8) formula}.  We have already demonstrated
that linear dependencies
exist between lattice sums (see \eqref{F(2,25) closed form}).  In certain cases it is
possible to relate new lattice sums to well-known examples.  Consider a forty-fifth degree modular equation due to Somos \cite[Entry $x_{45,4,12}$]{so}:
\begin{equation}\label{eta prod deg 45} 6 e_1 e_5 e_9 e_{45}=-e_1^2
e_5^2-2 e_3^2 e_{15}^2-9 e_{9}^2 e_{45}^2 + e_{3}^4 + 5 e_{15}^4.
\end{equation}
We were unable to prove \eqref{eta prod deg 45} by elementary
methods. Integrating \eqref{eta prod deg 45} leads to a linear
dependency between three lattice sums. We have
\begin{equation}\label{linear relation for F(5,9)}
9F(5,9)=45F(1,1)-50F(1,5).
\end{equation}
Both $F(1,1)$ and $F(1,5)$ equal values of
hypergeometric functions \cite{RV}, \cite{RZ1}.  Since the
difficult task of evaluating $F(1,5)$ is accomplished in
\cite{RZ1}, we easily obtain the following theorem.
\begin{theorem} Recall that $n(\alpha)$ is defined in \eqref{Mahler measure
definitions}.  We have
\begin{equation}\label{F(5,9) final formula}
\frac{108}{5\pi^2}F(5,9)=8n\left(3\sqrt[3]{2}\right)-9n\left(2\sqrt[3]{4}\right).
\end{equation}
\end{theorem}

 There are various additional formulas which follow from Boyd's Mahler measure conjectures.  A proof of Boyd's conductor $30$ conjectures would lead to
closed forms for both $F(2,15)$ and $F(2,5/3)$.  To make this explicit we use two relations.  First consider a four term modular equation which Somos highlighted in \cite{so2}:
\begin{equation}\label{Somos's second smallest identity}
e_{1}e_{3}e_{5}e_{15}+2e_{2}e_{6}e_{10}e_{30}=e_{1}e_{2}e_{15}e_{30}+e_{3}e_{5}e_{6}e_{10}.
\end{equation}
Integrating \eqref{Somos's second smallest identity}, and then using the evaluation $F(3,5)=4\pi^2 m(1)/15$ from \cite{RZ2}, leads to
\begin{equation}\label{linear relation for F(2,15)}
\begin{split}
F(2,15)+4F\left(2,\frac{5}{3}\right)
&=\frac{8\pi^2}{5}m(1).
\end{split}
\end{equation}
Next we require an unproven relation.  Boyd conjectured\footnote{See Table 2 in \cite{Bo1}.  In our notation, Boyd's entries correspond to values of $g(2-k)$} that for a conductor $30$ elliptic curve
\begin{equation*}
L(E_{30},2)\stackrel{?}{=}\frac{2\pi^2}{15}g(3),
\end{equation*}
where $g(\alpha)$ is defined in \eqref{Mahler measure definitions}.  The modularity theorem guarantees that $L(E_{30},2)=L(f_{30},2)$, where $f_{30}(e^{2\pi i \tau})$ is a weight-two cusp form on $\Gamma_{0}(30)$.
Somos has calculated a basis for the $1$-dimensional space of cusp forms on $\Gamma_{0}(30)$, and consequently the cusp form associated with conductor $30$ elliptic curves is given by
\begin{equation*}
f_{30}(q)=\eta(q^3)\eta(q^5)\eta(q^6)\eta(q^{10})-\eta(q)\eta(q^2)\eta(q^{15})\eta(q^{30}).
\end{equation*}
Upon integrating this cusp form, Boyd's conjecture becomes
\begin{equation}\label{linear relation 2 for F(2,15)}
F\left(2,\frac{5}{3}\right)-\frac{1}{4}F(2,15)\stackrel{?}{=}\frac{2\pi^2}{15}g(3).
\end{equation}
Combining \eqref{linear relation for F(2,15)} and \eqref{linear relation 2 for F(2,15)} leads to a pair of conjectural evaluations.
\begin{conjecture} Recall that $m(\alpha)$ and $g(\alpha)$ are defined in \eqref{Mahler measure definitions}. The following formulas are numerically true:
\begin{align}
\frac{15}{4\pi^2}F(2,15)&\stackrel{?}{=}3m(1)-g(3),\label{F(2,15) conjecture}\\
\frac{15}{\pi^2}F\left(2,\frac{5}{3}\right)&\stackrel{?}{=}3m(1)+g(3).\label{F(2,5/3) conjecture}
\end{align}
\end{conjecture}
Tracking backwards shows that a solution of either \eqref{F(2,15) conjecture} or \eqref{F(2,5/3) conjecture}
would settle Boyd's conductor $30$ Mahler measure conjectures.  Proofs remain out of reach, however we are optimistic that both identities
may eventually be established using Eisenstein series identities contained in \cite{BY}.

\section{Conclusion: Conductor $17$ elliptic curves}

A strong connection exists between lattice sums and Mahler measures, however
this relationship has limitations.  While our ultimate goal is to ``sum up" $F(b,c)$ for arbitrary values of $b$ and $c$, it is important to realize that
this would only settle a small portion of the conjectures in Boyd's paper \cite{Bo1}.
Conductor $17$ curves are the first cases in Cremona's list \cite{Cr}, where $L(E,2)$ probably does not reduce to values of $F(b,c)$.  If we let $E_{17}$ denote a conductor $17$ curve (we used $y^2+x y+y=x^3-x^2-x$), then
\begin{equation}\label{conductor 17}
\frac{17}{2\pi^2}L(E_{17},2)\stackrel{?}{=}m\left(\frac{(1+\sqrt{17})^2}{4}\right)-m\left(\sqrt{17}\right).
\end{equation}
We discovered \eqref{conductor 17} via numerical experiments involving elliptic dilogarithms.
The cusp form associated with conductor $17$ curves is stated in \cite{Fi}.  We have
\begin{equation}
f_{17}(q)=\frac{\eta(q)\eta^2(q^4)\eta^5(q^{34})}{\eta(q^2)\eta(q^{17})\eta^2(q^{68})}-\frac{\eta^5(q^2)\eta(q^{17})\eta^2(q^{68})}{\eta(q)\eta^2(q^4)\eta(q^{34})}.
\end{equation}
Since $L(E_{17},2)=L(f_{17},2)$, formula \eqref{conductor 17} can be changed into a complicated (!) elementary identity.  There does not seem to be an easy way to relate $L(E_{17},2)$ to Mahler measures of rational polynomials.  We surmise that this is the reason conductor $17$ curves never appear in Boyd's paper \cite{Bo1}.  Given the complexity of $f_{17}(q)$, we feel confident to conjecture that $L(E_{17},2)$ is linearly independent from values of $F(b,c)$ over $\mathbb{Q}$.

\begin{acknowledgements}
The authors thank David Boyd, Wadim Zudilin, and Bruce
Berndt for their useful comments and encouragement.  The authors are especially grateful to Michael Somos for the useful communications, and
for providing his list of modular equations.  Mat Rogers also thanks
the Max Planck Institute of Mathematics for their hospitality.
\end{acknowledgements}


\begin{thebibliography}{99}

\bibitem{BY}
\textsc{A.~Berkovich} and \textsc{H.~Yesilyurt},
Ramanujan's identities and representation of integers by certain binary and quaternary quadratic forms,
\emph{Ramanujan J.} \textbf{20} (2009), 375-408.

\bibitem{Be3}
\textsc{B.\,C.~Berndt},
\emph{Ramanujan's Notebooks, Part III}
(Springer-Verlag, New York, 1991).

\bibitem{Be4}
\textsc{B.\,C.~Berndt},
\emph{Ramanujan's Notebooks, Part IV}
(Springer-Verlag, New York, 1994).

\bibitem{Be5}
\textsc{B.\,C.~Berndt},
\emph{Ramanujan's Notebooks, Part V}
(Springer-Verlag, New York, 1998).



\bibitem{Bor}
\textsc{D.~Borwein}, \textsc{J. M.~Borwein} and \textsc{K.
F.~Taylor}, Convergence of lattice sums and Madelung's constant,
\emph{J. Math. Phys.} \textbf{26} (1985), no.~11, 2999--3009.

\bibitem{Bo1}
\textsc{D.\,W.~Boyd},
Mahler's measure and special values of $L$-functions,
\emph{Experiment. Math.} \textbf{7} (1998), 37--82.




\bibitem{Cr}
\textsc{J.\,E.~Cremona},
{}Algorithms for modular elliptic curves,
available at \texttt{http://www.warwick.ac.uk/~masgaj/ftp/data/}

\bibitem{De}
\textsc{C.~Deninger},
Deligne periods of mixed motives, $K$-theory and the entropy of certain $\mathbb Z^n$-actions,
\emph{J. Amer. Math. Soc.} \textbf{10} (1997), no.~2, 259--281.

\bibitem{Fi}
\textsc{S.~Finch},
Primitive cusp forms,
Preprint (2009).


\bibitem{GK}
\textsc{G.~K\"ohler},
\emph{Eta products and theta series identities}
(Springer-Verlag, Heidelberg, 2011).


\bibitem{KO}
\textsc{N.~Kurokawa} and \textsc{H.~Ochiai},
Mahler measures via crystalization,
\emph{Comment. Math. Univ. St. Pauli} \textbf{54} (2005), 121--137.


\bibitem{LR}
\textsc{M.~N.~Lal\'{i}n} and \textsc{M.~D.~Rogers},
Functional equations for Mahler measures of genus-one curves,
\emph{ Algebra and Number
Theory}, \textbf{1} (2007), no.~1, 87--117.


\bibitem{Ono2}
\textsc{K.~Ono},
The Web of Modularity: Arithmetic of the Coefficients of Modular Forms and q-series,
{\em American Mathematical Society, Providence, RI,\/} 2004.

\bibitem{Ra}
\textsc{S.~Ramanujan},
Modular equations and approximations to
$\pi$, [Quart. J. Math. \textbf{45} (1914), 350-372].
\textit{Collected papers of Srinivasa Ramanujan, 23-29, AMS Chelsea
Publ., Providence, RI, 2000.}

\bibitem{RV}
\textsc{F.~Rodriguez-Villegas},
Modular Mahler measures I,
in: \emph{Topics in number theory} (University Park, PA, 1997),
Math. Appl. \textbf{467} (Kluwer Acad. Publ., Dordrecht, 1999), 17--48.

\bibitem{Rgsubmit}
\textsc{M.~Rogers},
Hypergeometric formulas for lattice sums and Mahler measures,
\emph{Intern. Math. Res. Not.} (to appear),
preprint \texttt{arXiv:\,0806.3590 [math.NT]} (2008).

\bibitem{RZ1}
\textsc{M.~Rogers} and \textsc{W.~Zudilin},
{}From $L$-series of elliptic curves to Mahler measures,
preprint \texttt{arXiv:\,1012.3036 [math.NT]} (2010).

\bibitem{RZ2}
\textsc{M.~Rogers} and \textsc{W.~Zudilin},
{}On the mahler measure of $1+X+X^{-1}+Y+Y^{-1}$,
preprint \texttt{arXiv:\,1012.3036 [math.NT]} (2010).



\bibitem{so2}
\textsc{M.~Somos},
A Remarkable eta-product Identity, Preprint
(2008).

\bibitem{so}
\textsc{M.~Somos},
Dedekind eta function product identities,
available at \texttt{http://eta.math.georgetown.edu/}.

\bibitem{We} \textsc{E.~W.~Weisstein}, ``Ramanujan Constant." From MathWorld--A
Wolfram Web Resource.
\texttt{http://mathworld.wolfram.com/RamanujanConstant.html}

\bibitem{Zu} \textsc{W.~Zudilin},
{}Arithmetic hypergeometric series,
\emph{Russian Math. Surveys} \textbf{66} (2011), no.~1, 369–-420.


\end{thebibliography}
\end{document}